# SIMULATION AND APPROXIMATION OF LÉVY-DRIVEN STOCHASTIC DIFFERENTIAL EQUATIONS

NICOLAS FOURNIER

ABSTRACT. We consider the problem of the simulation of Lévy-driven stochastic differential equations. It is generally impossible to simulate the increments of a Lévy-process. Thus in addition to an Euler scheme, we have to simulate approximately these increments. We use a method in which the large jumps are simulated exactly, while the small jumps are approximated by Gaussian variables. Using some recent results of Rio [12] about the central limit theorem, in the spirit of the famous paper by Komlós-Major-Tsunády [10], we derive an estimate for the strong error of this numerical scheme. This error remains reasonable when the Lévy measure is very singular near 0, which is not the case when neglecting the small jumps.

In the same spirit, we study the problem of the approximation of a Lévy-driven S.D.E. by a Brownian S.D.E. when the Lévy process has no large jumps.

## 1. INTRODUCTION

Let $(Z_t)_{t\geq 0}$ be a one-dimensional square integrable Lévy process. Then for some $a \in \mathbb{R}$, some $b \in \mathbb{R}_+$, and some measure $\nu$ on $\mathbb{R}_*$ satisfying $\int_{\mathbb{R}_*} z^2 \nu(dz) < \infty$,

$$(1) \qquad Z_t = at + bB_t + \int_0^t \int_{\mathbb{R}_*} z \tilde{N}(ds, dz),$$

where $(B_t)_{t\geq 0}$ is a standard Brownian motion, independent of a Poisson measure $N(ds, dz)$ on $[0, \infty) \times \mathbb{R}_*$ with intensity measure $ds\nu(dz)$, and where $\tilde{N}$ is its compensated Poisson measure, see Jacod-Shiryaev [9].

We consider, for some $x \in \mathbb{R}$ and some function $\sigma : \mathbb{R} \mapsto \mathbb{R}$, the S.D.E.

$$(2) \qquad X_t = x + \int_0^t \sigma(X_{s-}) dZ_s.$$

Using some classical results (see e.g. Ikeda-Watanabe [4]), there is strong existence and uniqueness for (2) as soon as $\sigma$ is Lipschitz continuous: for any given couple $(B, N)$, there exists an unique càdlàg adapted solution $(X_t)_{t\geq 0}$ to (2). By *adapted*, we mean adapted to the filtration $(\mathcal{F}_t)_{t\geq 0}$ generated by $(B, N)$.

We consider two related problems in this paper. The first one, exposed in the next section, deals with the numerical approximation of the solution $(X_t)_{t\geq 0}$. The second one concerns the approximation of $(X_t)_{t\geq 0}$ by the solution to a Brownian S.D.E., when $Z$ has only very small jumps, and is presented in Section 3.

Our results are based on a recent work of Rio [12] that concerns the rate of convergence in the central limit theorem, when using the quadratic Wasserstein distance. This result, and its application to Lévy processes, is exposed in Section 4.

## 2. Numerical simulation

The first goal of this paper is to study a numerical scheme to solve (2). The first idea is to perform an Euler scheme $(X^n_{i/n})_{n \geq 0}$ with time-step $1/n$, see Jacod [5], Jacod-Protter [8], Protter-Talay [11] for rates of convergence. However, this is generally not a good scheme in practise, unless one knows how to simulate the increments of the underlying Lévy process, which is the case e.g. when $Z$ is a stable process.

We assume here that the Lévy measure $\nu$ is known explicitly: one can thus simulate random variables with law $\nu(dz)\mathbf{1}_A(z)/\nu(A)$, for any $A$ such that $\nu(A) < \infty$.

The first idea is to approximate the increments of $Z$ by $\widehat{\Delta}^{n,\epsilon}_i = Z^\epsilon_{i/n} - Z^\epsilon_{(i-1)/n}$, where $Z^\epsilon_t$ is the same Lévy process as $Z$ without its (compensated) jumps smaller than $\epsilon$. However, Asmussen-Rosinski [1] have shown that for a Lévy process with many small jumps, it is more convenient to approximate small jumps by some Gaussian variables than to neglect them. We thus introduce $\Delta^{n,\epsilon}_i = \widehat{\Delta}^{n,\epsilon}_i + U^{n,\epsilon}_i$, where $U^{n,\epsilon}_i$ is Gaussian with same mean and variance as the neglected jumps. The arguments of [1] concern only Lévy processes, and it does not seem so easy to apply such an idea to the simulation of SDEs.

Let us write $(\widehat{X}^{n,\epsilon}_{[nt]/n})_{t \geq 0}$ (resp. $(X^{n,\epsilon}_{[nt]/n})_{t \geq 0}$) for the Euler scheme using the approximate increments $(\widehat{\Delta}^{n,\epsilon}_i)_{i \geq 1}$ (resp. $(\Delta^{n,\epsilon}_i)_{i \geq 1}$). They of course have a similar computational cost.

Jacod-Kurtz-Méléard-Protter [7] have computed systematically the *weak* error for the *approximate* Euler scheme. In particular, they prove some very fine estimates of $\mathbb{E}[g(X^{n,\epsilon}_{[nt]/n})] - \mathbb{E}[g(X_t)]$ for $g$ smooth enough. The obtained rate of convergence is very satisfying.

Assume now that the goal is to approximate some functional of the path of the solution (e.g. $\sup_{[0,T]} |X_t|$). Then we have to estimate the error between the laws of the paths of the processes (not only between the laws of the time marginals). A common way to perform such an analysis is to introduce a suitable coupling between the numerical scheme $(X^{n,\epsilon}_{[nt]/n})_{t \geq 0}$ and the true solution $(X_t)_{t \geq 0}$, and to estimate the (discretized) *strong* error $\mathbb{E}[\sup_{[0,T]} |X^{n,\epsilon}_{[nt]/n} - X_{[nt]/n}|^2]$.

We refer to Jacod-Jakubowski-Mémin [6] for the rate of convergence of the discretized process $(X_{[nt]/n})_{t \geq 0}$ to the whole process $(X_t)_{t \geq 0}$.

Rubenthaler [13] has studied the strong error when neglecting small jumps. He obtains roughly $\mathbb{E}[\sup_{[0,T]} |\widehat{X}^{n,\epsilon}_{[nt]/n} - X_{[nt]/n}|^2] \simeq C_T(n^{-1} + \int_{|z| \leq \epsilon} z^2 \nu(dz))$ (if $b \neq 0$). For $\nu$ very singular near 0, the obtained precision is very low.

Our aim here is to study the strong error when using $X^{n,\epsilon}_{[nt]/n}$. We will see that the precision is much higher (see Subsection 2.4 below).

The main difficulty is to find a suitable coupling between the true increments $(Z_{i/n} - Z_{(i-1)/n})_{i \geq 1}$ and the approximate increments $(\Delta^{n,\epsilon}_i)_{i \geq 1}$: clearly, one considers $Z$, then one erases its jumps smaller than $\epsilon$, but how to build the additionnal Gaussian variable in such a way that it is a.s. close to the erased jumps? We will use a recent result of Rio [12], which gives some very precise rate of convergence for the standard central limit theorem in Wasserstein distance, in the spirit of Komlós-Major-Tsunády [10].

2.1. **Notation.** We introduce, for $\epsilon \in (0,1)$, $k \in \mathbb{N}$,

$$(3) \quad F_\epsilon(\nu) = \int_{|z|>\epsilon} \nu(dz), \quad m_k(\nu) = \int_{\mathbb{R}_*} |z|^k \nu(dz),$$

$$m_{k,\epsilon}(\nu) = \int_{|z|\leq\epsilon} |z|^k \nu(dz), \quad \delta_\epsilon(\nu) = \frac{m_{4,\epsilon}(\nu)}{m_{2,\epsilon}(\nu)}.$$

Observe that we always have $\delta_\epsilon(\nu) \leq \epsilon^2$ and $F_\epsilon(\nu) \leq \epsilon^{-2} m_2(\nu)$.

For $n \in \mathbb{N}$ and $t \geq 0$, we set $\rho_n(t) = [nt]/n$, where $[x]$ is the integer part of $x$.

2.2. **Numerical scheme.** Let $n \in \mathbb{N}$ and $\epsilon \in (0,1)$ be fixed. We introduce an i.i.d. sequence $(\Delta_i^{n,\epsilon})_{i\geq 1}$ of random variables, with

$$(4) \quad \Delta_1^{n,\epsilon} = a_{n,\epsilon} + b_{n,\epsilon} G + \sum_1^{N_{n,\epsilon}} Y_i^\epsilon,$$

where $a_{n,\epsilon} = (a - \int_{|z|>\epsilon} z\nu(dz))/n$, where $b_{n,\epsilon}^2 = (b^2 + m_{2,\epsilon}(\nu))/n$, where $G$ is Gaussian with mean 0 and variance 1, where $N_{n,\epsilon}$ is Poisson distributed with mean $F_\epsilon(\nu)/n$, and where $Y_1^\epsilon, Y_2^\epsilon, ...$ are i.i.d. with law $\nu(dz)\mathbf{1}_{|z|>\epsilon}/F_\epsilon(\nu)$. All these random variables are assumed to be independent. Then we introduce the scheme

$$(5) \quad X_0^{n,\epsilon} = x, \quad X_{(i+1)/n}^{n,\epsilon} = X_{i/n}^{n,\epsilon} + \sigma(X_{i/n}^{n,\epsilon})\Delta_{i+1}^{n,\epsilon} \quad (i \geq 0).$$

Observe that
- the cost of simulation of $\Delta_1^{n,\epsilon}$ is of order $1 + \mathbb{E}[N_{n,\epsilon}] = 1 + F_\epsilon(\nu)/n$, whence that of $(X_{\rho_n(t)}^{n,\epsilon})_{t\in[0,T]}$ is of order $Tn(1 + F_\epsilon(\nu)/n) = T(n + F_\epsilon(\nu))$, as in [13];
- $\Delta_{i+1}^{n,\epsilon}$ has the same law as $Z_{(i+1)/n}^\epsilon - Z_{i/n}^\epsilon + U_{n,\epsilon}$, where $U_{n,\epsilon}$ is Gaussian with same mean and variance as $\int_{i/n}^{(i+1)/n} \int_{|z|\leq\epsilon} z\tilde{N}(ds,dz)$ and where $Z_t^\epsilon = at + bB_t + \int_0^t \int_{|z|>\epsilon} z\tilde{N}(ds,dz)$.

2.3. **Main result.** We may now state our main result.

**Theorem 1.** *Assume that $\sigma : \mathbb{R} \mapsto \mathbb{R}$ is bounded and Lipschitz continuous. Let $\epsilon \in (0,1)$ and $n \in \mathbb{N}$. There is a coupling between a solution $(X_t)_{t\geq 0}$ to (2) and an approximated solution $(X_{\rho_n(t)}^{n,\epsilon})_{t\geq 0}$ as in Subsection 2.2 such that for all $T > 0$,*

$$\mathbb{E}\left[\sup_{[0,T]} |X_{\rho_n(t)} - X_{\rho_n(t)}^{n,\epsilon}|^2\right] \leq C_T \left(n^{-1} + n\delta_\epsilon(\nu)\right),$$

*where the constant $C_T$ depends only on $T, \sigma, a, b, m_2(\nu)$.*

The first bound $n^{-1}$ is due to the time discretization (Euler scheme), and the second bound $n\delta_\epsilon$ is due to the approximation of the increments of the Lévy process. As noted by Jacod [5], the first bound may be improved if there is no brownian motion $b = 0$ (but we have to work with some weaker norm).

We assume here that $m_2(\nu) < \infty$ for simplicity: this allows us to work in $L^2$. However, we believe that Theorem 1 allows one to show that in the general case where $\int_{\mathbb{R}_*} \min(z^2, 1)\nu(dz) < \infty$, the family $(n^{-1} + n\delta_\epsilon(\nu))^{-1} \sup_{[0,T]} |X_{\rho_n(t)} - X^{n,\epsilon}_{\rho_n(t)}|^2$ is tight: decompose the Lévy process $Z_t = Z_t^1 + Z_t^2$, where $Z^1$ satisfies our assumptions, and $Z^2$ is a compound Poisson process. Apply Theorem 1 between the jumps of $Z^2$, and paste the pieces... this might be complicated to write, but the principle is very simple.

2.4. **Optimisation.** Choose $\epsilon = 1/n$. Then recalling that $\delta_\epsilon(\nu) \leq \epsilon^2$, we get

$$\mathbb{E}\left[\sup_{[0,T]} |X_{\rho_n(t)} - X^{n,1/n}_{\rho_n(t)}|^2\right] \leq C_T/n,$$

for a mean cost to simulate $(X^{n,1/n}_{\rho_n(t)})_{t \in [0,T]}$ of order $T(n + F_{1/n}(\nu))$.
- We always have $F_\epsilon(\nu) \leq C\epsilon^{-2}$, so that the cost is always smaller than $Tn^2$.
- If $\nu(dz) \stackrel{z \to 0}{\simeq} |z|^{-1-\alpha} dz$ for some $\alpha \in (0,2)$, $F_\epsilon(\nu) \simeq \epsilon^{-\alpha}$, so that the cost is of order $Tn^{\max(1,\alpha)}$.

Still assume that $\nu(dz) \stackrel{z \to 0}{\simeq} |z|^{-1-\alpha} dz$, for some $\alpha \in (0,2)$. When neglecting the small jumps, the mean cost to get a mean squared error of order $1/n$ is of order $Tn^{\max(1,\alpha/(2-\alpha))}$ (see [13]), which is huge when $\alpha$ is close to 2. We observe that the present method is more precise as soon as $\alpha > 1$.

2.5. **Discussion.** The computational cost to get a given precision does not explode when the Lévy measure becomes very singular near 0. The more $\nu$ is singular at 0, the more there are jumps greater than $\epsilon$, which costs many simulations. But the more it is singular, the more the jumps smaller than $\epsilon$ are well-approximated by Gaussian random variables. These two phenomena are in competition, and we prove that the second one compensates (partly) the first one.

Our result involves a suitable coupling between the solution $(X_t)_{t \geq 0}$ and its approximation $(X_t^{n,\epsilon})_{t \geq 0}$. Of course, this is not very interesting in practise, since by definition, $(X_t)_{t \geq 0}$ is completely unknown. This is just an artificial way to estimate the rate of convergence *in law*, using a Wasserstein type distance.

The simulation algorithm can easily be adapted to the case of dimension $d \geq 2$. We believe that the result still holds. However, the result of Rio [12] is not known in the multidimensional setting (although it is believed to hold). We could use instead the results of Einmahl [2]. This would be much more technical, and would lead to a lower rate of convergence.

3. BROWNIAN APPROXIMATION

Consider the Lévy process introduced in (1), consider $x \in \mathbb{R}$, $\sigma : \mathbb{R} \mapsto \mathbb{R}$ Lipschitz continuous, and the unique solution $(X_t)_{t \geq 0}$ to (2). Recall (3), consider a Brownian motion $(W_t)_{t \geq 0}$, set

(6) $$\tilde{Z}_t = at + \sqrt{b^2 + m_2(\nu)} W_t,$$

which has the same mean and variance as $Z_t$. Let $(\tilde{X}_t)_{t \geq 0}$ be the unique solution to

(7) $$\tilde{X}_t = x + \int_0^t \sigma(\tilde{X}_{s-}) d\tilde{Z}_s.$$

**Theorem 2.** *Assume that $\sigma$ is Lipshitz continuous and bounded. Then it is possible to couple the solutions $(X_t)_{t\geq 0}$ to (2) and $(\tilde{X}_t)_{t\geq 0}$ to (7) in such a way that for all $p \geq 4$, all $T > 0$, all $n \geq 1$, recall (3)*

$$\mathbb{E}\left[\sup_{[0,T]} |X_t - \tilde{X}_t|^2\right] \leq C_{T,p}\left(n^{2/p-1} + m_p(\nu)^{2/p} + nm_4(\nu)\right),$$

*where $C_{T,p}$ depends only on $p, T, \sigma, a, b, m_2(\nu)$.*

If we only know that $m_4(\nu) < \infty$, then we choose $n = [m_4(\nu)^{-2/3}]$, and we get, at least when $m_4(\nu) \leq 1$, $\mathbb{E}\left[\sup_{[0,T]} |X_t - \tilde{X}_t|^2\right] \leq C_T m_4(\nu)^{1/3}$.

Consider a sequence of Lévy processes $(Z_t^\epsilon)_{t\geq 0}$ with drift $a$, diffusion coefficient $b$ and Lévy measure $\nu_\epsilon$, such that $z^2 \nu_\epsilon(dz)$ tends weakly to $\delta_0$. Then $\lim_{\epsilon \to 0} m_2(\nu_\epsilon) = 1$, while in almost all cases, $\lim_{\epsilon \to 0} m_p(\nu_\epsilon) = 0$ for some (are all) $p > 2$.
Consider the solution to $X_t^\epsilon = x + \int_0^t \sigma(X_{s-}^\epsilon) dZ_s^\epsilon$. Then it is well-known and easy to show that $(X_t^\epsilon)_{t\geq 0}$ tends in law to the solution of a Brownian S.D.E. Theorem 2 allows one to obtain a rate of convergence (for some Wasserstein distance). For example, we will immediately deduce the following corollary.

**Corollary 3.** *Assume that $\sigma$ is Lipshitz continuous and bounded. Assume that $\nu(\{|z| > \epsilon\}) = 0$ for some $\epsilon \in (0, 1]$. Then it is possible to couple the solutions $(X_t)_{t\geq 0}$ to (2) and $(\tilde{X}_t)_{t\geq 0}$ to (7) in such a way that for all $\eta \in (0, 1)$, all $T > 0$,*

$$\mathbb{E}\left[\sup_{[0,T]} |X_t - \tilde{X}_t|^2\right] \leq C_{T,\eta} \epsilon^{1-\eta}$$

*where $C_{T,\eta}$ depends only on $\eta, T, \sigma, a, b, m_2(\nu)$.*

The original motivation of this work was to estimate the error when approximating the Boltzmann equation by the Landau equation. The Boltzmann equation is a P.D.E. that can be related to a Poisson-driven S.D.E. (see Tanaka [14]), while the Landau equation can be related to a Brownian S.D.E. (see Guérin [3]). In the grazing collision limit, the S.D.E. related to the Boltzmann equation has only very small jumps. However, many additionnal difficulties arise for those equations. Furthermore, we are able to prove our results only in dimension 1, while the kinetic Boltzmann and Landau equations involve 3-dimensional S.D.E.s

## 4. Coupling results

Consider two laws $P, Q$ on $\mathbb{R}$ with finite variance. The Wasserstein distance $\mathcal{W}_2$ is defined by

$$\mathcal{W}_2^2(P, Q) = \inf \left\{\mathbb{E}\left[|X - Y|^2\right], \; \mathcal{L}(X) = P, \mathcal{L}(Y) = Q\right\}.$$

With an abuse of notation, we also write $\mathcal{W}_2(X, Y) = \mathcal{W}_2(X, Q) = \mathcal{W}_2(P, Q)$ if $\mathcal{L}(X) = P$ and $\mathcal{L}(Y) = Q$. We recall the following result of Rio [12, Theorem 4.1].

**Theorem 4.** *There is an universal constant $C$ such that for any sequence of i.i.d. random variables $(Y_i)_{i\geq 1}$ with mean $0$ and variance $\theta^2$, for any $n \geq 1$,*

$$\mathcal{W}_2^2\left(\frac{1}{\sqrt{n}} \sum_1^n Y_i, \mathcal{N}(0, \theta^2)\right) \leq C \frac{\mathbb{E}[Y_1^4]}{n\theta^2}$$

Here $\mathcal{N}(0,\theta^2)$ is the Gaussian distribution with mean 0 and variance $\theta^2$. Recall now (3).

**Corollary 5.** *Consider a pure jump centered Lévy process $(Y_t)_{t\geq 0}$ with Lévy measure $\mu$. In other words $Y_t = \int_0^t \int_{\mathbb{R}_*} z\tilde{M}(ds,dz)$, where $\tilde{M}$ is a compensated Poisson measure with intensity $ds\mu(dz)$. There is an universal constant $C$ such that*

$$\forall\, t \geq 0, \quad \mathcal{W}_2^2\left(Y_t, \mathcal{N}(0, tm_2(\mu))\right) \leq C\frac{m_4(\mu)}{m_2(\mu)}.$$

*Proof.* Let $t > 0$. For $n \geq 1$, $i \geq 1$, write $Y_i^n = n^{1/2}\int_{(i-1)t/n}^{it/n}\int_{\mathbb{R}_*} z\tilde{M}(ds,dz)$, whence $Y_t = n^{-1/2}\sum_1^n Y_i^n$. The $Y_i^n$ are i.i.d., centered, $\mathbb{E}[(Y_1^n)^2] = tm_2(\mu)$, and

$$\begin{aligned}
\mathbb{E}[(Y_1^n)^4] &= n^2\mathbb{E}\left[\left(\int_0^{t/n}\int_{\mathbb{R}_*} z^2 M(ds,dz)\right)^2\right] \\
&= n^2\mathbb{E}\left[\left(\int_0^{t/n}\int_{\mathbb{R}_*} z^2 \tilde{M}(ds,dz) + (t/n)m_2(\mu)\right)^2\right] \\
&= n^2\left[tm_4(\mu)/n + (tm_2(\mu)/n)^2\right] = ntm_4(\mu) + t^2 m_2(\mu).
\end{aligned}$$

Using Theorem 4, we get

$$\mathcal{W}_2^2\left(Y_t, \mathcal{N}(0, tm_2(\mu))\right) \leq C\frac{ntm_4(\mu) + t^2 m_2(\mu)}{ntm_2(\mu)} \xrightarrow{n\to\infty} C\frac{m_4(\mu)}{m_2(\mu)},$$

which concludes the proof. $\square$

This result is quite surprising at first glance: since the variances of the involved variables are $tm_2(\mu)$, it would be natural to get a bound that descreases to 0 as $t$ decreases to 0 (and that explodes for large $t$). Of course, we deduce the bound $\mathcal{W}_2^2(Y_t,\mathcal{N}(0,tm_2(\mu))) \leq C\min(m_4(\mu)/m_2(\mu), tm_2(\mu))$, but this is now optimal, as shown in the following example.

**Example.** Consider, for $\epsilon > 0$, $\mu_\epsilon = (2\epsilon^2)^{-1}(\delta_\epsilon + \delta_{-\epsilon})$, and the corresponding pure jump (centered) Lévy process $(Y_t^\epsilon)_{t\geq 0}$. It takes its values in $\epsilon\mathbb{Z}$. Observe that $m_2(\mu_\epsilon) = 1$ and $m_4(\mu_\epsilon) = \epsilon^2$. There is $c > 0$ such that for all $t \geq 0$, all $\epsilon > 0$, $\mathcal{W}_2^2(Y_t^\epsilon, \mathcal{N}(0,t)) \geq c\min(t,\epsilon^2) = c\min(m_4(\mu_\epsilon)/m_2(\mu_\epsilon), tm_2(\mu_\epsilon))$. Indeed,
• if $t \leq \epsilon^2$, then $\mathbb{P}(Y_t^\epsilon = 0) \geq e^{-t\mu_\epsilon(\mathbb{R}_*)} = e^{-t/\epsilon^2} \geq 1/e$, from which the lowerbound $\mathcal{W}_2^2(Y_t^\epsilon, \mathcal{N}(0,t)) \geq ct = c\min(t,\epsilon^2)$ is easily deduced;
• if $t \geq \epsilon^2$, use that $\mathcal{W}_2^2(Y_t^\epsilon, \mathcal{N}(0,t)) \geq \mathbb{E}[\min_{n\in\mathbb{Z}}|t^{1/2}G - n\epsilon|^2] = t\mathbb{E}[\min_{n\in\mathbb{Z}}|G - n\epsilon t^{-1/2}|^2]$, where $G$ is Gaussian with mean 0 and variance 1. Tedious computations show that there is $c > 0$ such that for any $a \in (0,1]$, $\mathbb{E}[\min_{n\in\mathbb{Z}}|G - na|^2] \geq (a/4)^2\mathbb{P}(G \in \cup_{n\in\mathbb{Z}}[(n+1/4)a, (n+3/4)a]) \geq ca^2$. Hence $\mathcal{W}_2^2(Y_t^\epsilon, \mathcal{N}(0,t)) \geq ct(\epsilon t^{-1/2})^2 = c\epsilon^2 = c\min(t,\epsilon^2)$.

## 5. Proof of Theorem 1

We recall elementary results about the Euler scheme for (2) in Subsection 5.1. We introduce our coupling in Subsection 5.2, which allows us to compare our scheme with the Euler scheme in Subsection 5.3. We conclude in Subsection 5.4. We assume in the whole section that $\sigma$ is bounded and Lipschitz continuous.

5.1. **Euler scheme.** We introduce the Euler scheme with step $1/n$ associated to (2). Let

(8) $$\Delta_i^n = Z_{i/n} - Z_{(i-1)/n} \quad (i \geq 1),$$

(9) $$X_0^n = x, \quad X_{(i+1)/n}^n = X_{i/n}^n + \sigma\left(X_{i/n}^n\right)\Delta_{i+1}^n \quad (i \geq 0).$$

The following result is classical.

**Proposition 6.** *Consider a Lévy process $(Z_t)_{t\geq 0}$ as in (1). For $(X_t)_{t\geq 0}$ the solution to (2) and for $(X_{i/n}^n)_{i\geq 0}$ defined in (8)-(9),*

$$\mathbb{E}\left[\sup_{[0,T]} |X_{\rho_n(t)} - X_{\rho_n(t)}^n|^2\right] \leq C_T/n,$$

*where $C_T$ depends only on $T, a, b, m_2(\nu)$, and $\sigma$.*

We sketch a proof for the sake of completeness.

*Proof.* Using the Doob and Cauchy-Scharz inequalities, we get, for $0 \leq s \leq t \leq T$,

(10) $$\mathbb{E}\left[\sup_{[s,t]} |X_u - X_s|^2\right] \leq C\mathbb{E}\Big[\left(a\int_s^t |\sigma(X_u)|du\right)^2 + \sup_{[s,t]}\left(b\int_s^v \sigma(X_u)dB_u\right)^2$$
$$+ \sup_{[s,t]}\left(\int_s^v \int_{\mathbb{R}_*} \sigma(X_{u-})z\tilde{N}(ds,du)\right)^2\Big]$$
$$\leq C_T \int_s^t (a^2 + b^2 + m_2(\nu))\|\sigma\|_\infty^2 du \leq C_T(t-s).$$

Observe now that $X_{\rho_n(t)}^n = x + \int_0^{\rho_n(t)} \sigma(X_{\rho_n(s)-}^n)dZ_s$. Setting $A_t^n = \sup_{[0,t]} |X_{\rho_n(s)} - X_{\rho_n(s)}^n|^2$, we thus get $A_t^n = \sup_{[0,t]} |\int_0^{\rho_n(s)} (\sigma(X_{u-}) - \sigma(X_{\rho_n(u)-}^n))dZ_u|^2$. Using the same arguments as in (10), then the Lipschitz property of $\sigma$ and (10), we get

$$\mathbb{E}[A_t^n] \leq C_T \int_0^{\rho_n(t)} (a^2 + b^2 + m_2(\nu))\mathbb{E}[(\sigma(X_s) - \sigma(X_{\rho_n(s)}^n))^2]ds$$
$$\leq C_T \int_0^t \mathbb{E}[(X_s - X_{\rho_n(s)})^2 + (X_{\rho_n(s)} - X_{\rho_n(s)}^n)^2]ds$$
$$\leq C_T \int_0^t (|s - \rho_n(s)| + \mathbb{E}[A_s^n])\, ds.$$

We conclude using that $|s - \rho_n(s)| \leq 1/n$ and the Gronwall Lemma. $\square$

5.2. **Coupling.** We now introduce a suitable coupling between the Euler scheme (see Subsection 5.1) and our numerical scheme (see Subsection 2.2). Recall (3).

**Lemma 7.** *Let $n \in \mathbb{N}$ and $\epsilon > 0$. It is possible to build two coupled families of i.i.d. random variables $(\Delta_i^n)_{i\geq 1}$ and $(\Delta_i^{n,\epsilon})_{i\geq 1}$, distributed respectively as in (8) and (4) in such a way that for each $i \geq 1$,*

$$\mathbb{E}[(\Delta_i^n - \Delta_i^{n,\epsilon})^2] \leq C\delta_\epsilon(\nu),$$

*where $C$ is an universal constant. Furthermore, for all $\epsilon > 0$, all $n \in \mathbb{N}$, all $i \geq 1$,*

$$\mathbb{E}[\Delta_i^n] = \mathbb{E}[\Delta_i^{n,\epsilon}] = \frac{a}{n}, \quad \mathrm{Var}[\Delta_i^n] = \mathrm{Var}[\Delta_i^{n,\epsilon}] = \frac{b^2 + m_2(\nu)}{n}.$$

*Proof.* It of course suffices to build $(\Delta_1^n, \Delta_1^{n,\epsilon})$, and then to take independent copies. Consider a Poisson measure $N(ds, dz)$ with intensity measure $ds\nu(dz)\mathbf{1}_{\{|z| \leq \epsilon\}}$ on $[0, \infty) \times \{|z| \leq \epsilon\}$. Observe that $\int_0^t \int_{|z| \leq \epsilon} z\tilde{N}(ds, dz)$ is a centered pure jump Lévy process with Lévy measure $\nu_\epsilon(dz) = \mathbf{1}_{|z| \leq \epsilon}\nu(dz)$. Then we use Corollary 5 and enlarge the underlying probability space if necessary: there is a Gaussian random variable $G_1^{n,\epsilon}$ with mean 0 and variance $m_2(\nu_\epsilon)/n = m_{2,\epsilon}(\nu)/n$ such that $\mathbb{E}\left[|\int_0^{1/n} \int_{|z| \leq \epsilon} z\tilde{N}(ds, dz) - G_1^{n,\epsilon}|^2\right] \leq Cm_4(\nu_\epsilon)/m_2(\nu_\epsilon) = C\delta_\epsilon(\nu)$.

We consider a Brownian motion $(B_t)_{t \geq 0}$ and a Poisson measure $N$ with intensity measure $ds\nu(dz)\mathbf{1}_{\{|z| > \epsilon\}}$ on $[0, \infty) \times \{|z| > \epsilon\}$, independent of the couple $(G_1^{n,\epsilon}, \int_0^{1/n} \int_{|z| \leq \epsilon} z\tilde{N}(ds, dz))$ and we set

- $\Delta_1^n := a/n + bB_{1/n} + \int_0^{1/n} \int_{|z| \leq \epsilon} z\tilde{N}(ds, dz) + \int_0^{1/n} \int_{|z| > \epsilon} z\tilde{N}(ds, dz)$,
- $\Delta_1^{n,\epsilon} := a/n + bB_{1/n} + G_1^{n,\epsilon} + \int_0^{1/n} \int_{|z| > \epsilon} z\tilde{N}(ds, dz)$.

Then $\Delta_1^n$ has obviously the same law as $Z_{1/n} - Z_0$ (see (1) and (8)), while $\Delta_1^{n,\epsilon}$ has also the desired law (see (4)). Indeed, $bB_{1/n} + G_1^{n,\epsilon}$ has a centered Gaussian law with variance $b^2/n + m_{2,\epsilon}(\nu)/n = b_{n,\epsilon}^2$, and $a/n + \int_0^{1/n} \int_{|z| > \epsilon} z\tilde{N}(ds, dz) = a_{n,\epsilon} + \int_0^{1/n} \int_{|z| > \epsilon} zN(ds, dz)$. This last integral can be represented as in (4). Finally $\mathbb{E}[(\Delta_1^n - \Delta_1^{n,\epsilon})^2] \leq \mathbb{E}\left[|\int_0^{1/n} \int_{|z| \leq \epsilon} z\tilde{N}(ds, dz) - G_1^{n,\epsilon}|^2\right] \leq C\delta_\epsilon(\nu)$, and the mean and variance estimates are obvious. □

### 5.3. Estimates.

We now compare our scheme with the Euler scheme. To this end, we introduce some notation. First, we consider the sequence $(\Delta_i^n, \Delta_i^{n,\epsilon})_{i \geq 1}$ introduced in Lemma 7. Then we consider $(X_{i/n}^n)_{i \geq 0}$ and $(X_{i/n}^{n,\epsilon})_{i \geq 0}$ defined in (9) and (5). We introduce the filtration $\mathcal{F}_i^{n,\epsilon} = \sigma(\Delta_k^n, \Delta_k^{n,\epsilon}, k \leq i)$, and the processes (with $V_0^{n,\epsilon} = 0$)

$$Y_i^{n,\epsilon} = X_{i/n}^n - X_{i/n}^{n,\epsilon}, \quad V_i^{n,\epsilon} = \frac{a}{n}\sum_{k=0}^{i-1}[\sigma(X_{k/n}^n) - \sigma(X_{k/n}^{n,\epsilon})], \quad M_i^{n,\epsilon} = Y_i^{n,\epsilon} - V_i^{n,\epsilon}.$$

**Lemma 8.** *There is a constant $C$, depending only on $\sigma, a, b, m_2(\nu)$ such that for all $N \geq 1$,*

$$\mathbb{E}\left[\sup_{i=0,\ldots,N} |Y_i^{n,\epsilon}|^2\right] \leq Cn\delta_\epsilon(\nu)(1 + C/n)^N(1 + N^2/n^2).$$

*Proof.* We divide the proof into four steps.

*Step 1.* We prove that for all $i \geq 0$, $\mathbb{E}\left[|Y_i^{n,\epsilon}|^2\right] \leq Cn\delta_\epsilon(\nu)(1 + C/n)^i$. First,

$$\mathbb{E}[|Y_{i+1}^{n,\epsilon}|^2] = \mathbb{E}[|Y_i^{n,\epsilon}|^2] + \mathbb{E}[(\sigma(X_{i/n}^n)\Delta_{i+1}^n - \sigma(X_{i/n}^{n,\epsilon})\Delta_{i+1}^{n,\epsilon})^2]$$
$$+ 2\mathbb{E}\left[Y_i^{n,\epsilon}(\sigma(X_{i/n}^n)\Delta_{i+1}^n - \sigma(X_{i/n}^{n,\epsilon})\Delta_{i+1}^{n,\epsilon})\right] = \mathbb{E}[|Y_i^{n,\epsilon}|^2] + I_i^{n,\epsilon} + J_i^{n,\epsilon}.$$

Now, using Lemma 7 and that $(\Delta_{i+1}^n, \Delta_{i+1}^{n,\epsilon})$ is independent of $\mathcal{F}_i^{n,\epsilon}$, we deduce that

$$J_i^{n,\epsilon} = \frac{2a}{n}\mathbb{E}\left[Y_i^{n,\epsilon}(\sigma(X_{i/n}^n) - \sigma(X_{i/n}^{n,\epsilon}))\right] \leq \frac{C}{n}E[|Y_i^{n,\epsilon}|^2],$$

since $\sigma$ is Lipschitz continuous. Using now the Lipschitz continuity and the boundedness of $\sigma$, together with Lemma 7 and the independence of $(\Delta_{i+1}^n, \Delta_{i+1}^{n,\epsilon})$ with respect to $\mathcal{F}_i^{n,\epsilon}$, we get

$$I_i^{n,\epsilon} \leq C\mathbb{E}[|Y_i^{n,\epsilon}|^2 (\Delta_{i+1}^{n,\epsilon})^2] + C\mathbb{E}[(\Delta_{i+1}^{n,\epsilon} - \Delta_{i+1}^n)^2] \leq \frac{C}{n} E[|Y_i^{n,\epsilon}|^2] + C\delta_\epsilon(\nu).$$

Finally, we get

$$\mathbb{E}[|Y_{i+1}^{n,\epsilon}|^2] \leq (1+C/n) E[|Y_i^{n,\epsilon}|^2] + C\delta_\epsilon(\nu).$$

Since $Y_0^{n,\epsilon} = 0$, this entails that $E[|Y_i^{n,\epsilon}|^2] \leq C\delta_\epsilon(\nu)[1 + (1+C/n) + ... + (1+C/n)^{i-1}] \leq Cn\delta_\epsilon(\nu)(1+C/n)^i$.

*Step 2.* We check that for $N \geq 1$, $\mathbb{E}[\sup_{0,...,N} |V_i^{n,\epsilon}|^2] \leq Cn\delta_\epsilon(\nu)(1+C/n)^N N^2/n^2$. It suffices to use the Lipschitz property of $\sigma$, the Cauchy-Schwarz inequality, and then Step 1:

$$\mathbb{E}\left[\sup_{1,...,N} |V_i^{n,\epsilon}|^2\right] \leq C\mathbb{E}\left[\left(\frac{1}{n} \sum_0^{N-1} |Y_i^{n,\epsilon}|\right)^2\right] \leq C\frac{N}{n^2} \sum_0^{N-1} \mathbb{E}[|Y_i^{n,\epsilon}|^2]$$

$$\leq C\frac{N^2}{n^2} n\delta_\epsilon(\nu)(1+C/n)^N$$

*Step 3.* We now verify that $(M_i^{n,\epsilon})_{i\geq 0}$ is a $(\mathcal{F}_i^{n,\epsilon})_{i\geq 0}$-martingale. We have $M_{i+1}^{n,\epsilon} - M_i^{n,\epsilon} = \sigma(X_{i/n}^n)[\Delta_{i+1}^n - a/n] - \sigma(X_{i/n}^{n,\epsilon})[\Delta_{i+1}^{n,\epsilon} - a/n]$. The step is finished, since the variables $\Delta_{i+1}^n - a/n$ and $\Delta_{i+1}^{n,\epsilon} - a/n$ are centered and independent of $\mathcal{F}_i^{n,\epsilon}$.

*Step 4.* Using the Doob inequality and then Steps 1 and 2, we get

$$\mathbb{E}\left[\sup_{i=0,...,N} |M_i^{n,\epsilon}|^2\right] \leq C \sup_{i=0,...,N} \mathbb{E}\left[|M_i^{n,\epsilon}|^2\right]$$

$$\leq C \sup_{i=0,...,N} \mathbb{E}\left[|Y_i^{n,\epsilon}|^2\right] + C \sup_{i=0,...,N} \mathbb{E}\left[|V_i^{n,\epsilon}|^2\right]$$

$$\leq Cn\delta_\epsilon(\nu)(1+C/n)^N(1+N^2/n^2).$$

But now

$$\mathbb{E}\left[\sup_{i=0,...,N} |Y_i^{n,\epsilon}|^2\right] \leq C\mathbb{E}\left[\sup_{i=0,...,N} |M_i^{n,\epsilon}|^2\right] + C\mathbb{E}\left[\sup_{i=0,...,N} |V_i^{n,\epsilon}|^2\right],$$

which allows us to conclude. $\square$

Let us rewrite these estimates in terms of $X^n$ and $X^{n,\epsilon}$.

**Lemma 9.** *Consider the sequence $(\Delta_i^n, \Delta_i^{n,\epsilon})_{i\geq 1}$ introduced in Lemma 7, and then $(X_{i/n}^n)_{i\geq 0}$ and $(X_{i/n}^{n,\epsilon})_{i\geq 0}$ defined in (9) and (5). For all $T \geq 0$,*

$$\mathbb{E}\left[\sup_{[0,T]} |X_{\rho_n(t)}^n - X_{\rho_n(t)}^{n,\epsilon}|^2\right] \leq C_T n\delta_\epsilon(\nu),$$

*where $C_T$ depends only on $T, a, b, m_2(\nu), \sigma$.*

*Proof.* With the previous notation, $\sup_{[0,T]} |X_{\rho_n(t)}^n - X_{\rho_n(t)}^{n,\epsilon}| = \sup_{i=0,...,[nT]} |Y_i^{n,\epsilon}|$. Thus using Lemma 8, we get the bound $Cn\delta_\epsilon(\nu)(1+C/n)^{[nT]}(1+[nT]^2/n^2) \leq Cn\delta_\epsilon(\nu)e^{CT}(1+T^2)$, which ends the proof. $\square$

### 5.4. Conclusion. We finally give the

*Proof of Theorem 1.* Fix $n \in \mathbb{N}$ and $\epsilon > 0$. Denote by $Q(du, dv)$ the joint law of $(\Delta_1^n, \Delta_1^{n,\epsilon})$ built in Lemma 7, and write $Q(du, dv) = Q_1(du)R(u, dv)$, where $R(u, dv)$ is the law of $\Delta_1^{n,\epsilon}$ conditionnally to $\Delta_1^n = u$.

Consider a Lévy process $(Z_t)_{t \geq 0}$ as in (1), and $(X_t)_{t \geq 0}$ the corresponding solution to (2). Set, for $i \geq 0$, $\Delta_i^n = Z_{i/n} - Z_{(i-1)/n}$, and consider the Euler scheme $(X_{i/n}^n)_{i \geq 0}$ as in (9). For each $i \geq 1$, let $\Delta_i^{n,\epsilon}$ be distributed according to $R(\Delta_i^n, dv)$, in such a way that $(\Delta_i^{n,\epsilon})_{i \geq 1}$ is an i.i.d. sequence. Finally, let $(X_{i/n}^{n,\epsilon})_{i \geq 0}$ as in (5).

By this way, the processes $(X_t)_{t \geq 0}$, $(X_{i/n}^n)_{i \geq 0}$ and $(X_{i/n}^{n,\epsilon})_{i \geq 0}$ are coupled in such a way that we may apply Proposition 6 and Lemma 9. We get

$$\mathbb{E}\left[\sup_{[0,T]} |X_{\rho_n(t)} - X_{\rho_n(t)}^{n,\epsilon}|^2\right] \leq 2\mathbb{E}\left[\sup_{[0,T]} |X_{\rho_n(t)} - X_{\rho_n(t)}^n|^2 + \sup_{[0,T]} |X_{\rho_n(t)}^n - X_{\rho_n(t)}^{n,\epsilon}|^2\right]$$
$$\leq C_T[n^{-1} + n\delta_\epsilon(\nu)].$$

This concludes the proof. $\square$

### 6. Proofs of Theorem 2 and Corollary 3

We assume in the whole section that $\sigma$ is bounded and Lipschtiz continuous. We start with a technical lemma.

**Lemma 10.** *Let $(X_t)_{t \geq 0}$ and $(\tilde{X}_t)_{t \geq 0}$ be solutions to (2) and (7). Then for $p \geq 2$, for all $t_0 \geq 0$, all $h \in (0, 1]$,*

$$\mathbb{E}\left[\sup_{[t_0, t_0+h]} |X_t - X_{t_0}|^p\right] \leq C_p(h^{p/2} + hm_p(\nu)), \quad \mathbb{E}\left[\sup_{[t_0, t_0+h]} |\tilde{X}_t - \tilde{X}_{t_0}|^p\right] \leq C_p h^{p/2},$$

*where $C_p$ depends only on $p, \sigma, a, b, m_2(\nu)$.*

*Proof.* It clearly suffices to treat the case of $(X_t)_{t \geq 0}$. Let thus $p \geq 2$.
Using the Burholder-Davies-Gundy inequality and the boundedness of $\sigma$, we get

$$\mathbb{E}\left[\sup_{[t_0, t_0+h]} |X_t - X_{t_0}|^p\right] \leq C_p \mathbb{E}\left[\left(\int_{t_0}^{t_0+h} |a\sigma(X_s)|ds\right)^p\right]$$
$$+ C_p \mathbb{E}\left[\left(\int_{t_0}^{t_0+h} b^2\sigma^2(X_s)ds\right)^{p/2}\right] + C_p \mathbb{E}\left[\left(\int_{t_0}^{t_0+h} \int_{\mathbb{R}_*} \sigma^2(X_s)z^2 N(ds, dz)\right)^{p/2}\right]$$
$$\leq C_p h^p + C_p h^{p/2} + C_p \mathbb{E}\left[\left(\int_{t_0}^{t_0+h} \int_{\mathbb{R}_*} z^2 N(ds, dz)\right)^{p/2}\right] \leq C_p h^{p/2} + C_p \mathbb{E}[U_h^{p/2}],$$

where $U_t = \int_0^t \int_{\mathbb{R}_*} z^2 N(ds, dz)$. It remains to check that for $t \geq 0$, $\mathbb{E}[U_t^{p/2}] \leq C_p(t^{p/2} + t m_p(\nu))$. But, with $C_p$ depending on $m_2(\nu)$,

$$\mathbb{E}[U_t^{p/2}] = \int_0^t ds \int_{\mathbb{R}_*} \nu(dz) \mathbb{E}[(U_s + z^2)^{p/2} - U_s^{p/2}]$$

$$\leq C_p \int_0^t ds \int_{\mathbb{R}_*} \nu(dz) \mathbb{E}[z^2 U_s^{p/2-1} + |z|^p] \leq C_p \int_0^t \mathbb{E}[U_s^{p/2-1}] ds + C_p m_p(\nu) t$$

$$\leq C_p \int_0^t \mathbb{E}[U_s^{p/2}] \epsilon^{-1} ds + C_p(\epsilon^{p/2-1} + m_p(\nu)) t,$$

for any $\epsilon > 0$. Hence $\mathbb{E}[U_t^{p/2}] \leq C_p(\epsilon^{p/2-1} t + m_p(\nu) t) e^{C_p t/\epsilon}$ by the Gronwall Lemma. Choosing $\epsilon = t$, we conclude that $\mathbb{E}[U_t^{p/2}] \leq C_p(t^{p/2} + m_p(\nu) t)$. □

*Proof of Theorem 2.* We fix $n \geq 1$, $T > 0$, and $p \geq 4$.

*Step 1.* Using Lemma 5 (see also Lemma 7) we deduce that we may couple two i.i.d. families $(\Delta_i^n)_{i \geq 1}$ and $(\tilde{\Delta}_i^n)_{i \geq 1}$, in such a way that:
- $(\Delta_i^n)_{i \geq 1}$ has the same law as the increments $(Z_{i/n} - Z_{(i-1)/n})_{i \geq 1}$ of the Lévy process (1);
- $(\Delta_i^n)_{i \geq 1}$ has the same law as the increments $(\tilde{Z}_{i/n} - \tilde{Z}_{(i-1)/n})_{i \geq 1}$ of the Lévy process (6);
- for all $i \geq 1$, $\mathbb{E}[(\Delta_i^n - \tilde{\Delta}_i^n)^2] \leq C m_4(\nu)$ (we allow constants to depend on $m_2(\nu)$).

*Step 2.* We then set $X_0^n = \tilde{X}_0^n = x$, and for $i \geq 1$, $X_{i/n}^n = X_{(i-1)/n}^n + \sigma(X_{(i-1)/n}^n) \Delta_i^n$ and $\tilde{X}_{i/n}^n = \tilde{X}_{(i-1)/n}^n + \sigma(\tilde{X}_{(i-1)/n}^n) \tilde{\Delta}_i^n$. Using exactly the same arguments as in Lemmas 8 and 9, we deduce that $\mathbb{E}\left[\sup_{[0,T]} |X_{\rho_n(t)}^n - \tilde{X}_{\rho_n(t)}^n|^2\right] \leq C_T n m_4(\nu)$, where $C_T$ depends only on $T, \sigma, a, b, m_2(\nu)$.

*Step 3.* But $(X_{\rho_n(t)}^n)_{t \geq 0}$ is the Euler discretization of (2), while $(\tilde{X}_{\rho_n(t)}^n)_{t \geq 0}$ is the Euler discretization of (7). Hence using Step 2, Proposition 6 and a suitable coupling as in the final proof of Theorem 1, $\mathbb{E}\left[\sup_{[0,T]} |X_{\rho_n(t)} - \tilde{X}_{\rho_n(t)}|^2\right] \leq C_T(1/n + n m_4(\nu))$.

*Step 4.* We now prove that $\mathbb{E}\left[\sup_{[0,T]} |X_t - X_{\rho_n(t)}|^2\right] \leq C_{T,p}(n^{2/p-1} + m_p(\nu)^{2/p})$. We set $\Gamma_i = \sup_{[i/n, (i+1)/n]} |X_t - X_{\rho_n(t)}| = \sup_{[i/n, (i+1)/n]} |X_t - X_{i/n}|$. By Lemma 10, $\mathbb{E}[\Gamma_i^p] \leq C_p[(1/n)^{p/2} + m_p(\nu)/n]$. Thus, since $p \geq 2$,

$$\mathbb{E}\left[\sup_{[0,T]} |X_t - X_{\rho_n(t)}|^2\right] \leq \mathbb{E}\left[\sup_{1,\ldots,[nT]} \Gamma_i^2\right] \leq \mathbb{E}\left[\sup_{1,\ldots,[nT]} \Gamma_i^p\right]^{2/p} \leq \mathbb{E}\left[\sum_1^{[nT]} \Gamma_i^p\right]^{2/p}$$

$$\leq C_{T,p} n^{2/p} \left[(1/n)^{p/2} + m_p(\nu)/n\right]^{2/p},$$

which ends the step.

*Step 5.* Exactly as in Step 4, we get $\mathbb{E}\left[\sup_{[0,T]} |\tilde{X}_t - \tilde{X}_{\rho_n(t)}|^2\right] \leq C_{T,p} n^{2/p-1}$.

*Step 6.* Using Steps 3, 4 and 5, we deduce that with a suitable coupling, we have $\mathbb{E}[\sup_{[0,T]} |X_t - \tilde{X}_t|^2] \leq C_{T,p}(n^{2/p-1} + m_p(\nu)^{2/p} + n^{-1} + n m_4(\nu))$. □

We conclude the paper with the

*Proof of Corollary 3.* Since $\nu(\{|z| > \epsilon\}) = 0$, we deduce that $m_p(\nu) \leq m_2(\nu)\epsilon^{p-2}$, for any $p \geq 2$. Applying Theorem 2 and choosing $n = [\epsilon^{-p/(p-1)}]$, we get the bound

$$C_{T,p}\left(\epsilon^{(1-2/p)(p/(p-1))} + \epsilon^{(p-2)(2/p)} + \epsilon^{2-p/(p-1)}\right) \leq C_{T,p}(\epsilon^{1-1/(p-1)} + \epsilon^{2-4/p}).$$

Hence for $\eta \in (0,1)$, it is possible to get the bound $C_{T,\eta}\epsilon^{1-\eta}$, choosing $p$ large enough. □

**Acknowledgement.** I wish to thank Jean Jacod for fruitfull discussions.

## References

[1] S. Asmussen, J. Rosiński, *Approximations of small jumps of Lévy processes with a view towards simulation*, J. Appl. Probab. Volume 38, Number 2 (2001), 482-493.
[2] U. Einmahl, *Extensions of results of Komlos, Major, and Tusnady to the multivariate case*, J. Multivariate Anal. 28 (1989), no. 1, 20–68.
[3] H. Guérin, *Solving Landau equation for some soft potentials through a probabilistic approach*, Ann. Appl. Probab. 13 (2003), no. 2, 515–539.
[4] N. Ikeda, S. Watanabe, *Stochastic differential equations and diffusion processes,* Second edition. North-Holland Mathematical Library, 24. North-Holland Publishing Co., Amsterdam; Kodansha, Ltd., Tokyo, 1989.
[5] J. Jacod *The Euler scheme for Lévy driven stochastic differential equations: limit theorems*, Ann. Probab. 32 (2004), no. 3A, 1830–1872.
[6] J. Jacod, A. Jakubowski, J. Mémin *On asymptotic errors in discretization of processes*, Ann. Probab. 31 (2003), no. 2, 592–608.
[7] J. Jacod, T. Kurtz, S. Méléard, P. Protter *The approximate Euler method for Lévy driven stochastic differential equations*, Ann. Inst. H. Poincaré Probab. Statist. 41 (2005), 523-558.
[8] J. Jacod, P. Protter *Asymptotic error distributions for the Euler method for stochastic differential equations*, Ann. Probab. 26 (1998), no. 1, 267–307.
[9] J. Jacod, A.N. Shiryaev *Limit theorems for stochastic processes,* second edition, Springer-Verlag, Berlin, 2003.
[10] J. Komlós, P. Major, G. Tusnády, *An approximation of partial sums of independent rvs and the sample df I.*, Z. Wahrsch. verw. Gebiete 32 111-131 (1975).
[11] P. Protter, D. Talay *The Euler scheme for Lévy driven stochastic differential equations*, Ann. Probab. 25 (1997), no. 1, 393–423.
[12] E. Rio, *Upper bounds for minimal distances in the central limit theorem*, to appear in Ann. Inst. Henri Poincaré Probab. Stat., 2009.
[13] S. Rubenthaler, *Numerical simulation of the solution of a stochastic differential equation driven by a Lévy process*, Stochastic Process. Appl. 103 (2003), no. 2, 311–349.
[14] H. Tanaka, *Probabilistic treatment of the Boltzmann equation of Maxwellian molecules*, Z. Wahrsch. Verw. Gebiete 46 (1978/79), no. 1, 67–105.

Nicolas Fournier, LAMA UMR 8050, Faculté de Sciences et Technologies, Université Paris Est, 61 avenue du Général de Gaulle, 94010 Créteil Cedex, France, email: nicolas.fournier@univ-paris12.fr